\documentclass[12pt]{article}
\usepackage[french]{babel}
\usepackage[T1]{fontenc}
\usepackage{amsmath,amsfonts,amsthm,amssymb, latexsym,euscript,dsfont,graphicx,units,color}
\usepackage[latin1]{inputenc}
\usepackage[font=footnotesize]{subfig}
\usepackage{pdfsync}

\begin{document}


\begin{center}
{\Large
	{\sc  Modélisation de réseaux de régulation de gènes par processus déterministes par morceaux}
}
\bigskip

 Aurélie Muller-Gueudin $^{1}$
\bigskip

{\it
$^{1}$ 
Institut Elie Cartan Nancy, Nancy-Université\\
Boulevard des Aiguillettes B.P. 239 \\
F-54506 Vandoeuvre lès Nancy\\
\texttt{aurelie.muller@iecn.u-nancy.fr}
}
\end{center}
\bigskip


{\bf R\'esum\'e.} Pour représenter l'évolution des espèces moléculaires dans un réseau de gènes, le modèle le plus classique est le processus de Markov à sauts. Ce modèle a l'inconvénient d'être long à simuler en raison de la rapidité et du grand nombre de réactions chimiques. Nous proposons des modèles approximatifs, basés sur les processus déterministes par morceaux, permettant de raccourcir les temps de simulation. Dans un article récent, nous avons montré rigoureusement la convergence des premiers modèles (processus de Markov à sauts) vers les seconds (processus déterministes par morceaux). Dans l'exposé, nous n'entrerons pas dans les détails de cette justification rigoureuse, mais nous montrerons des exemples d'application à des réseaux de gènes très simples (modèle de Cook et modèle du phage Lambda). 
\smallskip

{\bf Mots-cl\'es.} Réseau de gènes, Processus déterministes par morceaux, Processus de Markov à sauts, Simulation. 
\bigskip\bigskip

{\bf Abstract.} The molecular evolution in a gene regulatory network is classically modeled by Markov jump processes. However, the direct simulation of such models is extremely time consuming. Indeed, even the simplest Markovian model, such as the production module of a single protein involves tens of variables and biochemical reactions and an equivalent number of parameters. We study the asymptotic behavior of multiscale stochastic gene networks using
weak limits of Markov jump processes. The results allow us to propose new models with reduced execution times. In a new article, we have shown that, depending on the time and concentration scales of
the system, the Markov jump processes could be approximated by piecewise deterministic processes. We give some applications of our results for simple gene networks (Cook's model and Lambda-phage model). 
\smallskip

{\bf Keywords.} Gene networks, Piecewise deterministic processes, Markov jump processes, Simulation.

\section{Présentation du problème }
En biologie moléculaire, les réseaux de gènes sont définis par un ensemble de réactions chimiques entre les espèces moléculaires présentes dans une cellule considérée. Il est maintenant établi que la dynamique de ces réseaux est stochastique : les systèmes de réactions chimiques sont modélisés par des processus de Markov à sauts homogènes.  
 
Soit un ensemble de réactions chimiques notées $R_r$, pour $r\in \mathcal R$; l'ensemble $\mathcal R$ est supposé fini. Ces réactions modifient les quantités d'espèces moléculaires présentes dans la cellule. Chaque espèce moléculaire est notée $i$, pour $i\in S=\{1,\dots,M\}$. Le nombre de molécules de l'espèce moléculaire $i$ est noté $n_i$ et l'état du système est décrit par le vecteur $X=(n_1,\ldots, n_M)\in\mathbb N^{M}$.
Chaque réaction $R_r$ change l'état du système de la manière suivante :  $X\mapsto X+\gamma_r$, avec $\gamma_r\in\mathbb Z^M$. Le vecteur $\gamma_r$ est le saut associé à la réaction $R_r$. La réaction $R_r$ a lieu avec un taux $\lambda_r(X)$ qui dépend de l'état du système. 

Cette évolution est décrite par un processus de Markov à sauts. Les instants de sauts, notés $(T_j)_{j\ge 1}$ vérifient $T_0=0$, $T_j= \tau_1+\dots +\tau_j$, où $(\tau_k)_{k\ge 1}$ est une suite de variables aléatoires indépendantes et telles que 
$$
\mathbb P(\tau_i>t)= \exp\Big(-\displaystyle\sum_{r\in {\mathcal R}} \lambda_r(X(T_{i-1}))t\Big).
$$
Au temps $T_i$, la réaction $r\in {\mathcal R}$ a lieu avec probabilité
$\displaystyle\frac{\displaystyle\lambda_r\left(X(T_{i-1})\right)}{\displaystyle\sum_{r\in {\mathcal R}} \lambda_r\left(X(T_{i-1})\right)}$ et l'état du système change selon l'équation $X\to X+\gamma_r$, c'est-à-dire:
$$
X(T_{i})= X(T_{i-1})+\gamma_{r}.
$$
Ce processus de Markov a pour générateur : $A f(X) = \sum_{r\in \mathcal R} \left[ f(X+\gamma_r)-f(X)\right]\lambda_r(X),$ pour des fonctions $f$ appartenant au domaine du générateur.

Ces modèles de Markov à sauts ont des temps d'exécution très longs. Par exemple, même un simple réseau de gènes contient des dizaines de variables et  de paramètres. Nous avons montré (Crudu \textit{et al.}~(2007); Radulescu \textit{et al.}~(2012)) que ces modèles pouvaient être approchés (via des convergence en loi (Billingsley~(1999)) par des processus markoviens déterministes par morceaux. Ces modèles approchés ont amélioré les temps d'exécution. 

\section{Nos résultats}
Dans les applications, les nombres de molécules sont de différentes échelles : certaines espèces sont en grand nombre, et d'autres en faible nombre. En conséquence, nous avons décomposé l'ensemble des espèces en deux sous-ensembles, notés $C$ et $D$ de cardinaux $M_C$ et $M_D$. De même, l'état du système est noté 
$X=(X_C,X_D)$, et les sauts notés $\gamma_r=(\gamma_r^C,\gamma_r^D)$. Pour  $i\in D$, $n_i$ est d'ordre $1$, tandis que pour $i\in C$, $n_i$ est proportionnel à $N$ où $N$ est grand. Nos résultats asymptotiques correspondent à $N\to +\infty$. Notons également $\displaystyle x_C=\frac1N X_C$ et $x=(x_C,X_D)$. De même, l'ensemble des réactions se décompose suivant les espèces mises en jeu dans les réactions : $\mathcal R=\mathcal R_D\cup \mathcal R_C\cup \mathcal R_{DC}$. Une réaction dans $\mathcal R_D$ (resp. $\mathcal R_C$) produit ou consomme uniquement des espèces rares, c'est-à-dire de $D$ (resp. des espèces fréquentes, c'est-à-dire de $C$). De même, les taux de réaction dans $\mathcal R_{D}$ (resp. $\mathcal R_{C}$) dépendent uniquement de $X_{D}$ (resp. $x_{C}$). Une réaction dans $\mathcal R_{DC}$ a un taux qui dépend à la fois de $x_C$ et $X_D$ et produit ou consomme des espèces à la fois rares ($D$) et fréquentes ($C$). 

Les taux des réactions $r\in\mathcal R_{C}$ sont également grands et d'ordre $N$, et nous posons $\tilde \lambda_{r}= \frac{\lambda_{r}}{N}$. Ceci signifie que la variable fréquente $x_C$ est fréquemment impliquée dans des réactions chimiques. Supposons pour l'instant que les réactions dans $\mathcal R_{D}$ ou $\mathcal R_{DC}$ ont un taux d'ordre 1.

Ces nouvelles variables obtenues par changement d'échelle obéissent encore à un processus de Markov à saut, dont le générateur est :
$$
\begin{array}{ll}
\tilde{\mathcal{A} }f(x_C,X_D) &=\displaystyle  \sum_{r\in \mathcal R_C} \left[ f(x_C+\frac1N \gamma_r^C,X_D)-f(x_C,X_D)
\right]N\tilde\lambda_r(x_C)\\
\\
&\displaystyle +\sum_{r\in \mathcal R_{DC}} \left[ f(x_{C}+\frac1N\gamma_r^C,X_{D}+\gamma_{r}^D)-f(x_{C},X_{D})\right]\lambda_r(x_{C},X_{D})\\
\\
&\displaystyle +\sum_{r\in \mathcal R_{D}} \left[ f(x_{C}, X_{D}+\gamma_r^D)-f(x_{C},X_{D})\right]\lambda_r(X_{D}).
\end{array}
$$
Si $N\to+\infty$, Kurtz~(1971,1978) a montré que les espèces rares et les espèces fréquentes se découplent, c'est-à-dire obéissent à des processus qui leur sont propres. En effet, les réactions dans $\mathcal R_{DC}$ ne sont pas suffisamment fréquentes pour changer le comportement de la variable rapide $x_{C}$. A la limite, la variable $x_C$ obéit à un système différentiel, qui fonctionne sans influence de la variable discrète $X_D$. La variable discrète $X_D$, quant à elle, obéit à un processus de sauts, indépendant de la variable continue $x_C$.

Notre travail a été de considérer des système plus généraux contenant d'autres types de réactions, et aboutissant à d'autres systèmes limites. 
Selon les différentes échelles de temps de réactions et de concentration des espèces, nous avons démontré rigoureusement quatre types de limites pour le processus de Markov à sauts : des processus continus et déterministes par morceaux (Davis~(1993)), des processus déterministes par morceaux avec des sauts sur la variable 'continue', des processus déterministes par morceaux moyennés, et des processus déterministes par morceaux avec des sauts singuliers sur la variable 'continue'.

\section{Illustration}
\subsection{Le phage $\lambda$}

Le phage $\lambda$ est un parasite de la bactérie E.Coli. L'état du phage $\lambda$ est donné par le vecteur $X = (C,C_2,D,D_1,D_2) $, où $C$ et $C_1$ représentent une protéine et son dimère, produits par le phage, et $D,D_1,D_2$ représentent des sites promoteurs sur l'ADN du phage, respectivement non occupés, en simple ou en double occupation. Nous avons $D+D_1+D_2=$constante. Le phage se développe au sein de la bactérie E. Coli selon les réactions chimiques suivantes:
$2C \underset{k_{-1}}{\overset{k_1}{\rightleftharpoons}} C_2$,
$D+C_2 \underset{k_{-2}}{\overset{k_2}{\rightleftharpoons}} D_1$,
$D_1+C_2 \underset{k_{-3}}{\overset{k_3}{\rightleftharpoons}} D_2$,
 $D_1 \xrightarrow{k_4} D_1  + nC$,  $C \xrightarrow{k_5}$. La dernière équation signifie qu'une protéine $C$ meurt. 

\begin{figure}[hpbf]
\begin{center}
\includegraphics[width=6cm]{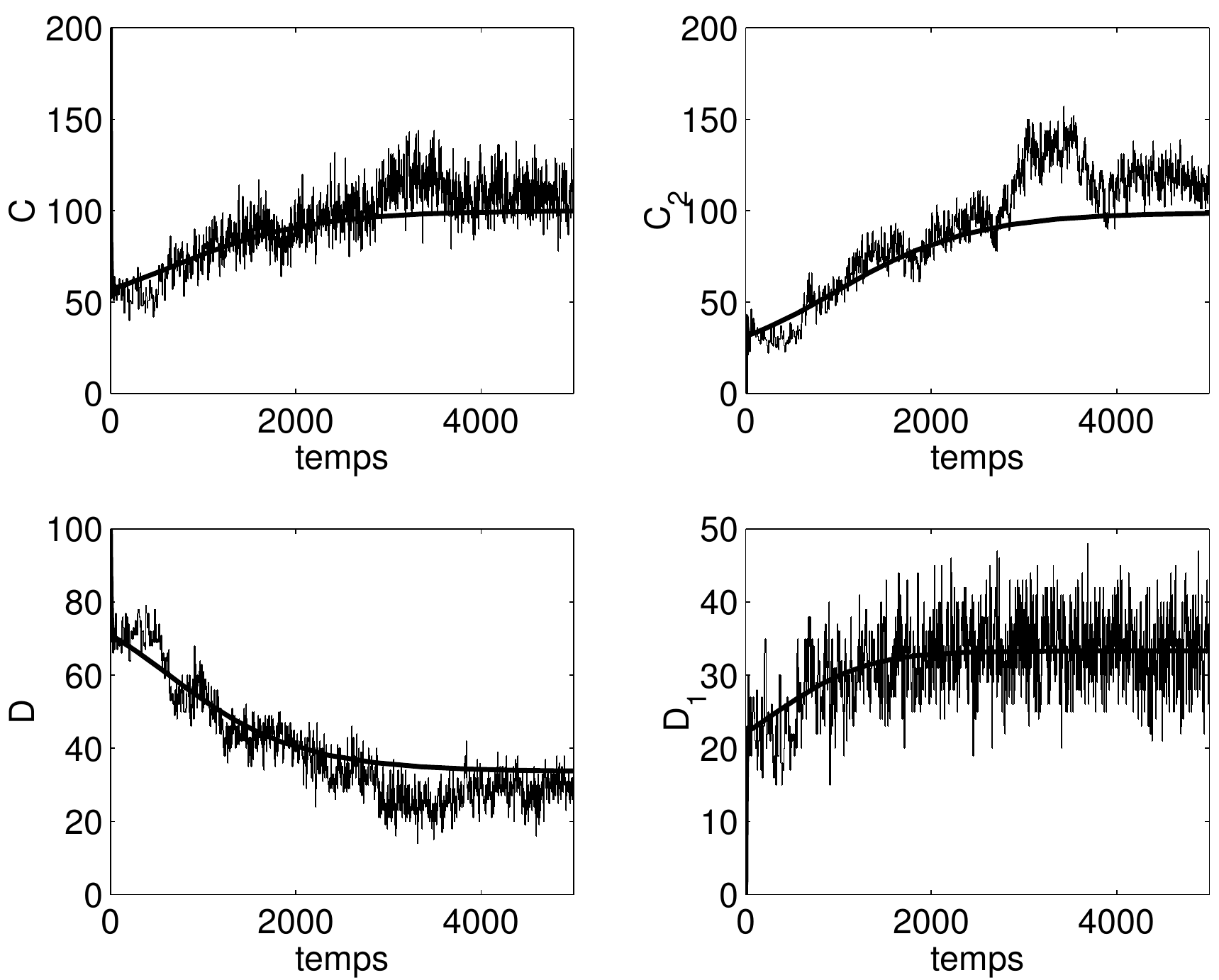}a)
\includegraphics[width=6cm]{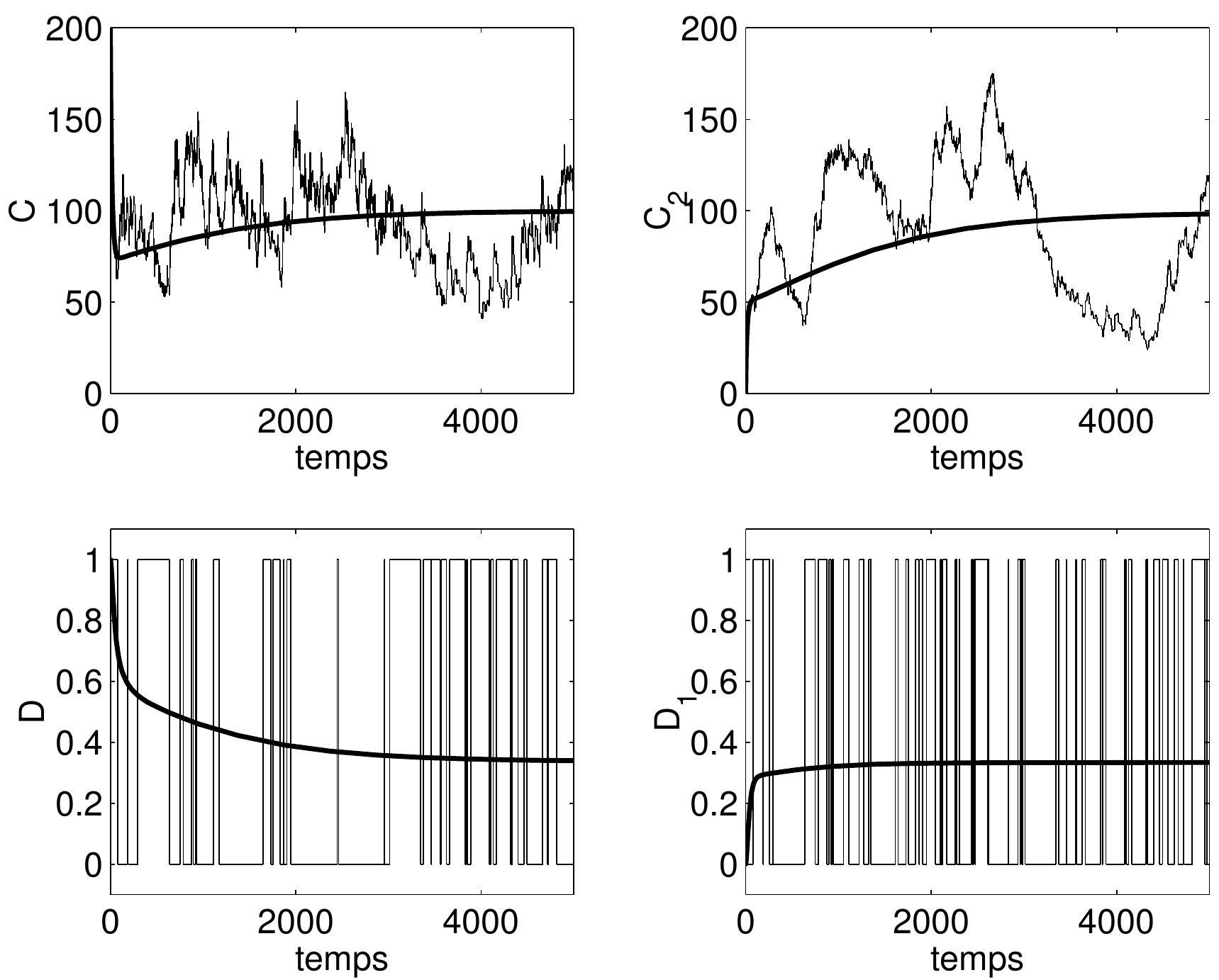}b)
\caption{Modèle de phage $\lambda$,  la limite déterministe est
marquée en trait continu. a) molécules et sites en grands nombres;
b) un seul site. Les paramètres sont a)
$k_1^{\pm}=k_2^{\pm}=k_3^{\pm}=0.1, k_4=0.006, k_5=0.01 $ b)
$k_1^{\pm}=k_3^{\pm}=0.01, k_4=0.3, k_5=0.005 $. }
\label{phage}
\end{center}
\end{figure}

Si les espèces moléculaires sont en grands nombres et ont des sauts rapides (ce qui se vérifie par certaines conditions sur les constantes cinétiques $k_i$), nous pouvons approcher la trajectoire de $X_t$ par une trajectoire déterministe. 
Les fluctuations autour des trajectoires déterministes sont illustrées dans la figure \ref{phage}~a).

En réalité, toutes les molécules ne sont pas en grand nombre, car $D,D_1,D_2$ valent soit 0 soit 1. Dans la figure ~\ref{phage}~b) nous avons simulé des trajectoires sous cette hypothèse. Les processus déterministes par morceaux sont des bonnes approximations pour les trajectoires de $C,C_2$. En effet, les molécules $C,C_2$ sont en grand nombre, on peut leur appliquer la limite déterministe entre deux sauts des variables discrètes. L'évolution de $(D,D_1)$ peut être décrite par un processus de Markov à sauts sur l'ensemble $\{0, 1 \}^2$. Lorsque le site est non occupé ou en double occupation $D=1,D_1=0$ ou
$D=0,D_1=0$, il n'y a pas de production de $C$ qui tend à
s'équilibrer avec son dimère. Lorsque le site est en simple
occupation $D=0,D_1=1$ il y a production de $C$. Le caractère déterministe par morceaux de la
dynamique des variables $C,C_2$ ne signifie pas que les fluctuations
de ces variables sont absentes. Bien au contraire, ces variables
sont soumises à des fluctuations importantes (Fig.~\ref{phage}~b).
Cette possibilité conduit a une conclusion biologique : la source
des fluctuations n'est pas nécessairement la petitesse du nombre moyen de
molécules observées mais pourrait être la petitesse du nombre de sites.

\subsection{Modèle de Cook}

Le modèle de Cook modélise les phénomènes de haploinsuffisance (invalidité de la moitié du nombre total de copies d'un gène). Il peut être décrit par le système de réactions suivant: $G \underset{k_{-1}}{\overset{k_1}{\rightleftharpoons}} G^*$, $G^*
\xrightarrow{k_2} G^* + P$, $P \xrightarrow{k_3}\cdot$

Dans le modèle de Cook $G,G^*$ sont des versions dormante et
active d'un gène, $P$ est la protéine traduite.

Ce système de réactions conserve la quantité $G+G^*=G_0$. Le régime
de haploinsuffisance est défini par une faible valeur $G_0$. Pour
simplifier considérons $G_0=1$, ce qui signifie une seule copie
valide du gène. Dans cette situation $G,G^* \in \{0,1 \}$. Plus
généralement, le modèle de Cook pourrait décrire d'autres cas
d'activité intermittente en biologie moléculaire. Par exemple, dans
le fonctionnement de voies de signalisation,  $G,G^*$ peuvent être
considérées comme des variables cachées à deux valeurs paramètrisant
la dynamique d'un système moléculaire dans deux situations (présence
et absence d'une molécule clé). Des extensions sont possibles à des
variables cachées à plusieurs valeurs discrètes (le comportement du
système pourrait dépendre de la présence et de l'absence de
certaines molécules).

Si les paramètres cinétiques $k_i$ vérifient certaines conditions, et si le nombre de
molécules $P$ est grand, on peut considérer que la dynamique de $P$
est déterministe pour une valeur de $G^*$ fixée. On arrive au modèle
déterministe par morceaux suivant:
$$
\frac{ dP}{dt} = -k_3 P + k_2 G^*(t)
$$

où $G^*(t)$ est un processus de Markov à espace d'états $E=\{0,1\}$
et fonction d'intensité:
$$
\lambda(G^*) = \left\{ \begin{array}{ll} k_1 & \text{si} \quad G^*=0
\\ k_{-1} & \text{si} \quad G^*=1 \end{array}\cdot \right.
$$

\begin{figure}[hpbf]
\begin{center}
\includegraphics[width=8cm]{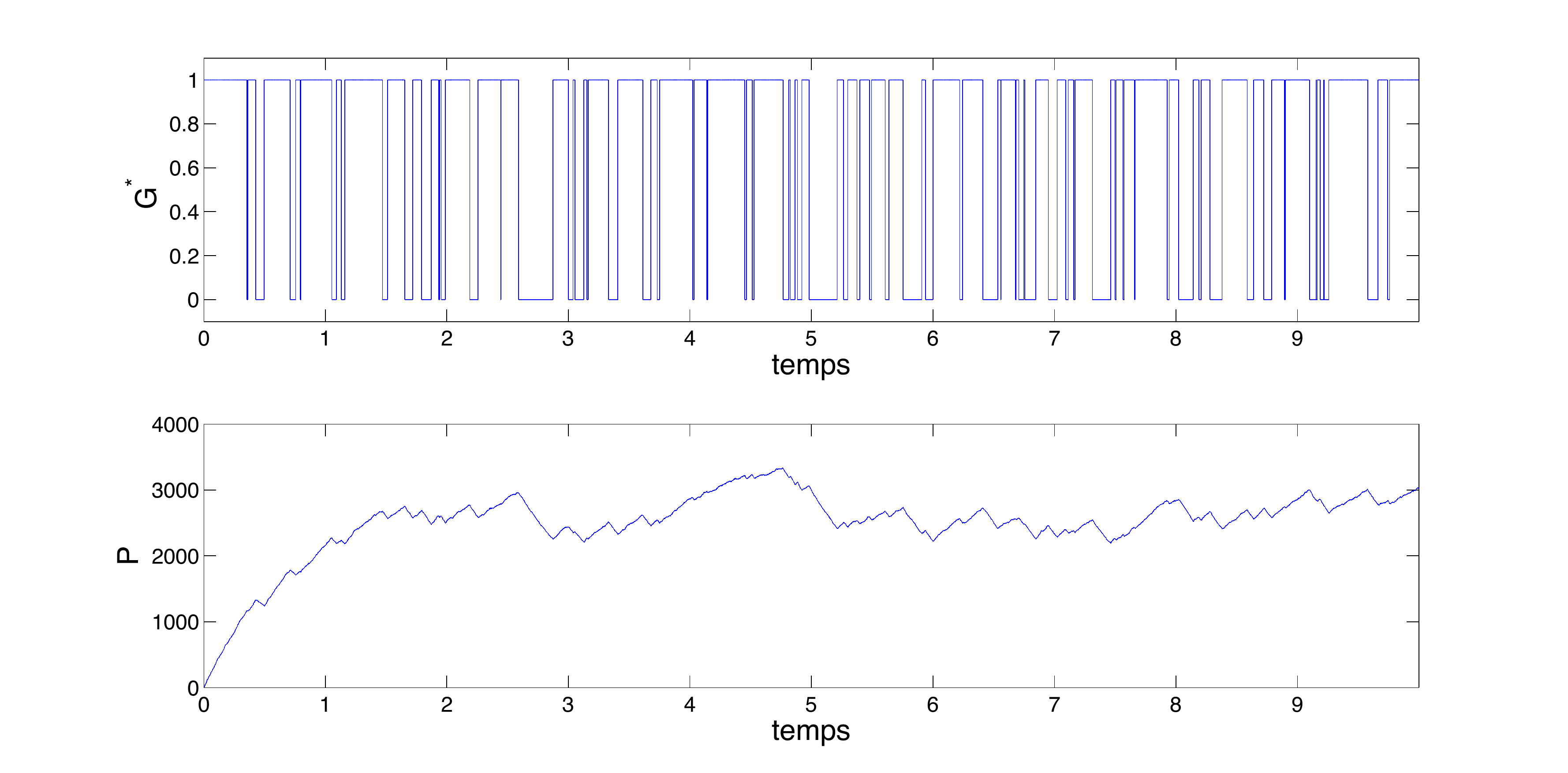}
\caption{Modèle de Cook. Les paramètres sont $k_1=20,k_{-1}=10, k_2=4000, k_3=1$.}
\label{Cook}
\end{center}
\end{figure}

Une trajectoire est simulée sur la figure \ref{Cook}. 

\section{Conclusion}
\begin{itemize}
\item 
Les simulations précédentes montrent que les temps d'exécutions des modèles déterministes par morceaux sont bien plus courts que ceux des modèles de Markov à sauts. En effet, il n'est plus nécessaire de simuler TOUTES les réactions chimiques, et donc TOUS les sauts du processus markovien. L'approximation déterministe par morceaux, justifiée théoriquement dans notre article, est une alternative économique aux modèles markoviens de sauts. 
\item 
Les trajectoires obtenues avec ces modèles approchés (les modèles déterministes par morceaux) sont similaires aux trajectoires obtenues en simulant toutes les réactions des modèles de Markov à sauts. 
\item Nous avons généralisé les résultats de Kurtz~(1971,1978). En effet, Kurtz a proposé des modèles déterministes qui ne sont pas adaptés dans les situations où les espèces moléculaires conservent un comportement stochastique. Nos modèles, déterministes par morceaux, ont permis de palier à cette lacune. Les exemples du phage $\lambda$ et de Cook illustrent ce point. 
\end{itemize}



\section*{Bibliographie}

\noindent [1] Billingsley, P. (1999), {\it Convergence of Probability Measures}, Wiley Series in Probability Statistics.

\noindent[2] Crudu, A., Debussche, A., Muller, A. et Radulescu, O. (2012),
Convergence of stochastic gene networks to hybrid piecewise deterministic processes, {\it A paraître dans Annals of Applied Probability}.

\noindent[3] Davis, M. (1993), {\it Markov Models and Optimization}, Chapman and Hall.

\noindent[4] Ethier, S. and Kurtz, T. (1986), {\it Markov processes. {C}haracterization and Convergence}, Wiley Series in Probability Statistics.

\noindent[5] Kurtz, T. (1971), Limits theorems for sequences of jump markov processes approximating ordinary differential processes, {\it J.Appl.Prob.}, 8,344--356.

\noindent[6] Kurtz, T. (1978), Strong approximation theorems for density dependent markov chains, {\it Stoch. Proc. Appl.}, 6,223--240.

\noindent[7] Radulescu, O., Muller, A. et Crudu, A. (2007), Théorèmes limites pour les processus de Markov à sauts. Synthèse de résultats et applications en biologie moléculaire, {\it Technique et Science Informatiques}, 3-4,441--467

\end{document}